\documentclass[12pt]{article}

\usepackage{amsmath,amssymb,amsfonts}

\newcommand{\Ann}{\mbox{Ann}}
\newcommand{\Hom}{\mbox{Hom}}
\newcommand{\Ext}{\mbox{Ext}}
\newcommand{\Tor}{\mbox{Tor}}
\newcommand{\Spec}{\mbox{Spec}\,}

\newcommand{\Supp}{\mbox{Supp}\,}
\newcommand{\gr}{\mbox{grade}\,}
\newcommand{\depth}{\mbox{depth}\,}
\renewcommand{\dim}{\mbox{dim}\,}
\newcommand{\cmd}{\mbox{cmd}\,}
\newcommand{\imp}{\mbox{imp}\,}

\newcommand{\pd}{\mbox{proj.dim}\,}

\newcommand{\gd}{\mbox{G--dim}\,}

\newcommand{\T}{\mathrm}

\newcommand{\lo}{\longrightarrow}

\newcommand{\para}{\paragraph}

\newcommand{\fm}{\frak{m}}
\newcommand{\fp}{\frak{p}}
\newcommand{\fq}{\frak{q}}

\begin{document}

\begin{center}

\LARGE{\bf Cohen--Macaulayness of Tensor Product}

\end{center}

\begin{center}

{\bf Leila Khatami and Siamak Yassemi\footnote{This research was
supported in part by a grant from IPM.}}
\end{center}

\begin{center}
{\it  Department of Mathematics, University of Tehran\\
and
\\ Institute for Studies in Theoretical Physics and Mathematics (IPM).}
\end{center}

\vspace{.3in}

\begin{abstract}

Let $(R,\fm)$ be a commutative Noetherian local ring. Suppose that
$M$ and $N$ are finitely generated modules over $R$ such that $M$
has finite projective dimension and such that $\Tor^R_i(M,N)=0$
for all $i>0$. The main result of this note gives a condition on
$M$ which is necessary and sufficient for the tensor product of
$M$ and $N$ to be a Cohen--Macaulay module over $R$, provided $N$
is itself a Cohen--Macaulay module.

\end{abstract}

\vspace{.2in}

{\bf 1991 Mathematics subject classification.} 13C14, 13D45, 13H10.

{\bf Key words and phrases.} Cohen--Macaulay modules.

\baselineskip=18pt

\vspace{.3in}

\section*{0. Introduction} Throughout this note $(R,\fm)$ is a
commutative Noetherian local ring with non-zero identity and the
maximal ideal $\fm$. By $M$ and $N$ we always mean non-zero
finitely generated $R$--modules. The projective dimension of
module $M$ is denoted by $\pd M$.

The well known notion ``grade of $M$'', grade $M$, has been
introduced by Rees, see [{\bf 7}], as the least integer $t\geq 0$
such that $\Ext_R^t(M,R)\neq 0$. In [{\bf 9}] we have defined the
``grade of $M$ and $N$'', grade $(M,N)$, as the least integer
$t\geq 0$ such that $\Ext_R^t(M,N)\neq 0$. One of the main
results of this note is:

Let $N$ be a Cohen--Macaulay $R$--module and let $M$ be an
$R$--module with finite projective dimension. If $\Tor_i^R(M,N)=0$
for all $i>0$ then $M\otimes_R N$ is Cohen--Macaulay if and only
if $\gr (M,N)=\pd M$.

This theorem can be considered as a generalization of the
well--known statement:

({\bf T1})\, Let $R$ be a Cohen--Macaulay local ring and let $M$
be a finite $R$--module with finite projective dimension. Then $M$
is Cohen--Macaulay if and only if $\gr  M=\pd M$.

On the other hand the following statement from Yoshida can be
concluded from our result:

Yoshida [{\bf 10}; Prop (2.4)] ``Suppose that $\gr  M=\pd
M(<\infty)$ and that $N$ is a maximal Cohen--Macaulay $R$--module
(that is depth $N=\dim N=\dim R$). Then $M\otimes_R N$ is
Cohen--Macaulay and $\dim M\otimes_R N=\dim M$.''

In another theorem of the first section we improve a theorem from
Kawasaki:

Kawasaki [{\bf 5}; Theorem 3.3(i)] ``Let $R$ be a Cohen--Macaulay
local ring and let $K$ be a canonical module of $R$. Let $M$ be  a
finite $R$--module of finite projective dimension. Then
$M\otimes_R K$ is Cohen--Macaulay if and only if $M$ is
Cohen--Macaulay.''

The following statement generalizes Kawasaki's theorem:

Let $R$ be a Cohen--Macaulay local ring and let $K$ be a canonical
module of $R$. If $M$ is an $R$--module with finite Gorenstein
dimension, then $M\otimes_R K$ is Cohen--Macaulay if and only if
$M$ is Cohen--Macaulay.

In the above statement the Gorenstein dimension is an invariant
for finite modules which was introduced by Auslander, in [{\bf
1}]. It is a  finer invariant than projective dimension in the
sense that $\gd M\leq \pd M$ for every finite non-zero $R$--module
$M$ and equality holds when $\pd  M <\infty$. There exist modules
with finite Gorenstein dimension which have infinite projective
dimension.

In the second section we consider Serre's condition. We  say $M$
satisfies Serre's condition $(S_n)$, for a non--negative integer
$n$, when for every $\fp \in \Supp M$ the following inequality
holds:
$$\depth M_{\fp}\geq \min (n,\dim M_{\fp}).$$
Obviously every Cohen--Macaulay module satisfies $(S_n)$ for all
non-negative integers $n$.

The main result of section 2 is:

Let $M$ and $N$ be $R$--modules such that $\Tor_i^R(M,N)=0$ for
all $i>0$. If projective dimension of $M$ is finite and
$M\otimes_R N$ satisfies $(S_n)$, then so does $N$. This result
generalizes [{\bf 10}; Prop.(4.1)].

\vspace{.3in}

\section*{1. Cohen--Macaulayness.}

\vspace{.2in}

\noindent{\bf Definition 1.1} We define
$$\gr  (M,N)=\inf \{i\, |\,\Ext_R^i
(M,N)\neq 0\}$$ Since $M$ is finite, using [{\bf 3}; 1.2.10]
we have that

$$\begin{array}{rl}
\gr  (M,N)&\, =\inf \{\depth N_{\fp}\, |\, {\fp}\in \Supp
M\}\\
&\, =\inf \{\depth N_{\fp}\, |\, {\fp}\in \Supp M\cap \Supp N\}.
\end{array}$$

The second equality holds because the depth of the zero module is
defined to be infinite.

\vspace{.2in}

\noindent{\bf Proposition 1.2} [{\bf 9}; Theorem 2.1] the
following inequalities hold:
\begin{verse}

(a) $\depth N-\dim M\leq\gr (M,N)$

(b) If $\Supp M\subseteq\Supp N$ then $\gr (M,N)\leq \dim N-\dim
M$
\end{verse}

\vspace{.2in}

For a finite $R$--module $M$ of finite projective dimension, the
invariant $\imp M$, imperfection of $M$, is defined to be $\pd M -
\gr  M$. This is, using Auslander--Buchsbaum equality, equal to
$\depth R-\depth M-\gr  M$.

\vspace{.2in}

\noindent{\bf Definition 1.3} For finite $R$--modules $M$ and $N$
(which may have infinite projective dimensions) we define
$\imp(M,N)=\depth N-\depth M-\gr  (M,N)$ (This may be negative).

It is clear that if $\pd  M <\infty$, then $\imp M=\imp(M,R)$.\\
By $\cmd M$ we mean the difference $\dim M-\depth M$.

\vspace{.2in}

\noindent{\bf Proposition 1.4} The following inequalities hold:
\begin{verse}

(a) $\imp(M,N)\leq \cmd M$

(b) If $\Supp M \subseteq \Supp N$, then $\cmd M\leq
\imp(M,N)+\cmd N$.

\end{verse}

\noindent{\em Proof}. This is clear from Proposition 1.2 and the
definition.

\para{Corollary 1.5} Let $N$ be a Cohen--Macaulay $R$--module and $\Supp
M\subseteq\Supp N$. Then $\cmd M=\imp(M,N)$; in particular the
module $M$ is a Cohen--Macaulay module if and only if
$\imp(M,N)=0$.\hfill$\square$

\vspace{.2in}

(T1) says that over a Cohen--Macaulay local ring $R$, the
$R$--module $M$ with finite projective dimension is
Cohen--Macaulay if $\Ext_R^i (M,R)=0$ for $i\neq\pd M$. The
following corollary is a generalization of (T1).

\vspace{.1in}

\noindent{\bf Corollary 1.6} Let $N$ be a Cohen--Macaulay
$R$--module with $\depth N=\depth R$. Let $M$ have finite
projective dimension and $\Supp M \subseteq \Supp N$. Then $M$ is
Cohen--Macaulay if and only if $\Ext_R^i (M,N)=0$ for $i\neq\pd
M$.

\vspace{.1in}

\noindent{\em Proof.} Note that $\pd  M =\sup\{i|
\Ext^i_R(M,N)\neq 0\, \, \mbox{for any}\, \, N\}$, cf. [{\bf 6}]
and so it is always greater than or equal to $\gr  (M,N)$.

$$\begin{array}{rl}
\imp(M,N) &\, =\depth N-\depth M-\gr (M,N)\\
&\, =\depth R-\depth M-\gr (M,N)\\
&\, =\pd  M -\gr (M,N)
\end{array}$$

Now the claim is clear from Corollary 1.5 .\hfill$\square$

\vspace{.2in}

Recall that a finite $R$--module $M$ with finite projective
dimension is called perfect if $\pd M=\gr  M$.

\vspace{.1in}

\noindent{\bf Definition 1.7} Let $M$ and $N$ be $R$--modules with
$\pd  M<\infty$. We say that $M$ is $N$-perfect if $\pd  M =\gr
(M,N)$.

\vspace{.2in}

In the proof of the following statements we use the well--known
result:

({\bf T2})\, Let $M$ and $N$ be finite $R$--modules with $\pd  M
<\infty$. If $\Tor_i^R(M,N)=0$ for all $i>0$, then we have the
equality $\depth M\otimes_R N=\depth N-\pd  M $.

\vspace{.1in}

\noindent{\bf Theorem 1.8} Let $N$ be a Cohen--Macaulay
$R$--module and let $M$ be an $R$--module with finite projective
dimension. If $\Tor_i^R (M,N)=0$ for all $i>0$, then $M\otimes_R
N$ is Cohen--Macaulay if and only if $M$ is $N$--perfect.

\vspace{.1in}

\noindent{\em Proof}. We claim that $M$ is $N$--perfect if and
only if $\imp(M\otimes_R N,N)=0$, and then the assertion will be
clear from Corollary 1.5 . We know that $\depth M\otimes_RN=\depth
N-\pd  M$. On the other hand

$$\begin{array}{rl}
\gr (M\otimes_R N,N) &\, =\inf \{\depth N_{\fp}\, |\, {\fp}\in\Supp M\otimes_RN\}\\
&\, =\inf \{\depth N_{\fp}\, |\, {p}\in\Supp M\cap \Supp N\}\\
&\, =\gr (M,N).
\end{array}$$

Then we have the equality $\imp(M\otimes_R N,N)=\pd  M -\gr
(M,N)$, which proves our claim.\hfill$\square$

\vspace{.2in}

Now [{\bf 10}; 2.4] can be deduced from the above theorem, for
when $N$ is maximal Cohen--Macaulay and $\pd  M <\infty$ by [{\bf
10}; 2.2] we have that $\Tor_i^R(M,N)=0$ for all $i>0$. For every
${\fp}\in \Supp N$, the $R_{\fp}$--module $N_{\fp}$ is maximal
Cohen--Macaulay module and, then $\depth N_{\fp}=\dim R_{\fp}\geq
\depth R_{\fp}$ and hence we have inequalities $$\gr  M\leq \gr
(M,N)\leq \pd  M. $$ This means that every perfect module is
$N$--perfect.

\vspace{.1in}

\noindent{\bf Definition 1.9} A finite $R$--module $N$ is said to
be of Gorenstein dimension zero, $\gd N=0$, if and only if

\begin{verse}

(a) $\Ext_R^i(N,R)=0$ for $i>0$.

(b) $\Ext_R^i (\Hom_R(N,R),R)=0$ for $i>0$.

(c) The canonical map $N\lo \Hom_R(\Hom_R (N,R),R)$ is an
isomorphism.

\end{verse}

For a non--negative integer $n$, the $R$--module $N$ is said to be
of Gorenstein dimension at most $n$, if and only if there exists
an exact sequence $$0\lo G_n\lo G_{n-1}\lo\ldots \lo G_0\lo N\lo 0
$$ where $\gd G_i=0$ for $0\leq i\leq n$. If such a sequence
does not exist then $\gd N=\infty$.

\vspace{.2in}

\noindent{\bf Lemma 1.10} [{\bf 2}; 3.7, 3.14, and 4.12] If $\gd
M<\infty$ then the following hold:

\begin{verse}

(a) $\gd M+\depth M=\depth R$.

(b) $\gd M=\sup\{t\, |\,\Ext_R^t(M,R)\neq 0\}.$

(c) $\Tor_i^R(M,P)=0$  for all $i>\gd M$ and all modules $P$ with
finite projective dimension.

\end{verse}

\vspace{.2in}

The following theorem improves Kawasaki's result [{\bf 5};
3.3(i)].

\vspace{.2in}

\noindent{\bf Theorem 1.11} Let $R$ be a Cohen--Macaulay local
ring and let $K$ be a canonical module of $R$. If $M$ is an
$R$--module with finite Gorenstein dimension, then $M\otimes_R N$
is Cohen--Macaulay if and only if $M$ is Cohen--Macaulay.

\vspace{.1in}

\noindent{\em Proof}. Proposition [{\bf 4}; 2.5] says that
$\Tor_i^R(M,K)=0$ for $i>0$, and then since injective dimension of
$K$ is finite we have that $\depth M\otimes_RK=\depth K-\gd M$,
cf. [{\bf 8}; 2.13]. Then $\imp (M\otimes_R K,K)=\gd K-\gr
(M\otimes_R K,K)$. Since $\Supp K=\Spec R$ we have that

$$\begin{array}{rl}
\gr (M\otimes_R K,K)&\, =\inf \{\depth
K_{\fp}\, |\, {\fp}\in\Supp M\otimes K\} \\
&\, = \inf \{\depth K_{\fp}\, |\, {\fp}\in \Supp M\}
\end{array}$$

But since $\depth K_{\fp}=\depth R_{\fp}$ for all ${\fp}\in\Supp
K=\Spec R$ we have that $\gr (M\otimes_R K,K)=\gr  M$. The claim
of the theorem is now clear from Corollary 1.5 and the fact that
over a Cohen--Macaulay local ring $R$, the $R$--module $M$ with
$\gd M<\infty$ is Cohen--Macaulay if and only if $\gr  M=\gd M$,
cf. [{\bf 9}].\hfill$\square$

\vspace{.3in}

\section*{2. Serre Conditions.} First recall that for a non--negative integer
$n$, we say that a finite $R$--module $M$ satisfies Serre's
condition $(S_n)$ if $\depth  M_{\fp}\geq \min (n,\dim M_{\fp})$
for every $\fp\in\Supp M$ or equivalently if $M_{\fp}$ is a
Cohen--Macaulay $R_{\fp}$-module for every $\fp\in\Supp M$ such
that $\depth   M_{\fp}<n$.

We also recall the intersection theorem:

({\bf T3})\, Let $M$ and $N$ be finite $R$--modules with $\pd  M
<\infty$. We have the inequality $\dim N\leq \pd  M +\dim
(M\otimes_R N)$.

\vspace{.1in}

\noindent{\bf Theorem 2.1} Let $N$ be a finite $R$--module which
satisfies $(S_n)$. Let $M$ be an $N$-perfect $R$--module with
$t=\pd  M\leq n$, such that $\Tor_i^R(M,N)=0$ for $i>0$. Then
$M\otimes_R N$ satisfies $(S_{n-t})$.

\vspace{.1in}

\noindent{\em Proof}. For every $\fp\in\Supp (M\otimes_R N)$ it is
clear that $$\gr (M,N)\leq \gr (M_{\fp}, N_{\fp})\leq\pd
M_{\fp}\leq \pd M =t.$$ Since $M$ is $N$--perfect, $M_{\fp}$ is
$N_{\fp}$--perfect with $\pd  M_{\fp}=t$. From Proposition 1.2 we
have that $t=\gr (M_{\fp}, N_{\fp})=\gr ((M\otimes_R N)_{\fp},
N_{\fp})\leq \dim N_{\fp}-\dim (M_R\otimes N)_{\fp}$.

On the other hand from the fact that $N$
satisfies $(S_n)$ we have the following inequality:

$$
\begin{array}{rl}
\depth  (M\otimes_R N)_{\fp}
&\, =\depth  (M_{\fp}\otimes_{R_{\fp}}N_{\fp})\\
&\, = \depth  N_{\fp}-\pd M _{\fp}\\
&\, = \depth  N_{\fp}-t\\
&\, = \min(n, \dim N_{\fp})-t\\
&\, =\min(n-t, \dim N_{\fp}-t)
\end{array}
$$

Now the assertion holds.\hfill$\square$

\vspace{.2in}

\noindent{\bf Corollary 2.2} If $R$ satisfies $(S_n)$, then every perfect
$R$--module with projective dimension $t$ (less than or equal to
$n$) satisfies $(S_{n-t})$.

\vspace{.2in}

It is well known that if a local ring admits a finite
Cohen--Macaulay module with finite projective dimension, then the
ring itself is Cohen--Macaulay.

In [{\bf 10}; 4.1] Yoshida has proved a more general statement, by
replacing ``being Cohen--Macaulay'' with ``satisfying Serre's
condition $(S_n)$''.

Our next two theorems improve those results by similar proofs. The
Theorem 2.3 is a special case of the Theorem 2.4 , and the proof
of it is only included because it is so simple.

\vspace{.1in}

\noindent{\bf Theorem  2.3} Let $M$ and $N$ be $R$--modules such
that $\Tor_i^R(M,N)=0$ for all $i>0$. If $\pd  M <\infty$ and
$M\otimes_R N$ is Cohen--Macaulay then so is $N$.

\vspace{.1in}

\noindent{\em Proof}. The intersection theorem (T3) gives the
inequality: $$\dim N\leq \dim M\otimes_R N+\pd M $$ On the other
hand (T2) gives the equality
$$\depth  N=\depth  M\otimes_R N+\pd M .$$ Since $\dim N\geq
\depth N$, the assertion is clear.\hfill$\square$

\vspace{.2in}

\noindent{\bf Theorem 2.4} Let $M$ and $N$ be $R$--modules such
that $\Tor_i^R (M,N)=0$ for all $i>0$. If $\pd M <\infty$ and
$M\otimes_R N$ satisfies $(S_n)$, then so does $N$.

\vspace{.1in}

\noindent{\em Proof}. Choose $\fp\in\Supp N$. There are two
cases:\\ The first case is when $\fp\in\Supp M$ and then
$\fp\in\Supp M\otimes_R N$.

If $\depth (M\otimes_R N)_{\fp}<n$ then $(M\otimes_R N)_{\fp}$ is
Cohen--Macaulay and by the Theorem 2.3 so is $N_{\fp}$.

If $\depth (M\otimes_R N)_{\fp} \geq n$, then $\depth N_{\fp} \geq
n$ because by (T2) $\depth N_{\fp}=\depth (M\otimes_R N)_{\fp}+\pd
M_{\fp}$.  The second case is when ${\fp}\not\in \Supp M$. Let
$\frak q$ be a minimal prime over the ideal $(\Ann M+\fp)$. From
(T3) we have the inequality
$$\dim{R_{\fq}}/{\fp R_{\fq}} \leq\pd M_{\fq}+\dim {M_{\fq}}/{\fp
M_{\fq}}=\pd  M _{\fq}.$$ Since $\fp R_{\fq}\in\Supp
{R_{\fq}}/{\fp R_{\fq}}$ we have that

$$\begin{array}{rl}
\depth N_{\fp} &\, \geq\gr ({R_{\fq}}/{\fp R_{\fq}}, N_{\fq})\\
&\, \geq \depth N_{\fq} -\dim {R_{\fq}}/{\fp R_{\fq}}\ (\T{proposition} 1.2)\\
&\, \geq \depth N_{\fq} -\pd M_{\fq}\\
&\, = \depth M_{\fq} \otimes_{R_{\fq}}N_{\fq}
\end{array}$$

If $\depth  M_{\fq} \otimes_{R_{\fq}}N_{\fq}<n$, then
$M_{\fq}\otimes_{R_{\fq}}N_{\fq}$ is Cohen--Macaulay and from
Theorem 2.3 we will have that $N_{\fq}$ is Cohen--Macaulay, then
so is $N_{\fp}\cong(N_{\fq})_{\fp R_{\fq}}$.

If $\depth  M_{\fq}\otimes_{R_{\fq}}N_{\fq}\geq n$ then the above
inequality guarantees that $\depth  N_{\fp}\geq n$.\hfill$\square$

\vspace{.3in}

\begin{center}

{\large\bf Acknowledgment}

\end{center}

The authors would like to thank the University of Tehran and the
University of Copenhagen for the facilities offered during the
preparation of this paper. Our thanks also to Professor Enochs,
University of Kentucky and Professor Foxby, University of
Copenhagen for their valuable comments on this paper.

\vspace{.3in}

\baselineskip=16pt

\begin{center}
{\bf References}
\end{center}

\begin{itemize}

\item[{[1]}] M.\ Auslander, {\em Anneaux de Gorenstein et torsion en
alg\`{e}bre commutative}, S\'{e}minaire d'alg\`{e}bre commutative
1966/67, notes by M.\ Mangeney, C.\ Peskine and L. Szpiro,
\'{E}cole Normale Sup\'{e}rieure de Jeunes Filles, Paris,
1967.\item[{[2]}] M. Auslander and M. Bridger, {\it Stable Module
Theory}, Memoris. Amer. Math. Soc., {\bf 94} (1969). \item[{[3]}]
W. Bruns and J. Herzog, {\it Cohen-Macaulay rings}, Cambridge
University Press, Cambridge (1993).
\item[{[4]}] H.-B. Foxby, {\it  Gorenstein modules and
related modules}, Math. Scand. {\bf 31} (1972),
267-285.
\item[{[5]}] T. Kawasaki, {\it Surjective-Buchsbaum
modules over Cohen-Macaulay local rings}, Math. Z., {\bf 218}
(1995), 191-205. \item[{[6]}] H. Matsumura, Commutative Ring
Theory, Cambridge University Press, Cambridge (1986). \item[{[7]}]
D. Rees, {\it The grade of an ideal or module}, Proc. Camb. Phil.
Soc., {\bf 53} (1957), 28-42.
\item[{[8]}] S. Yassemi, {\it G--dimension}, Math. Scand., {\bf 77}
(1995), 161--174.
\item[{[9]}]S. Yassemi, L. Khatami, and T. Sharif, {\it Grade and
Gorenstein dimension}, to appear in Comm. Algebra.
\item[{[10]}] K. Yoshida, {\it Tensor Products of perfect
modules and maximal surjective Buchsbaum modules}, J.
Pure. appl. Algebra, {\bf 123} (1998), 313-326.

\end{itemize}
\end{document}